\documentclass[12pt]{amsart} 
\usepackage{amssymb} 
\usepackage{colordvi} 
\usepackage{graphicx} 
\usepackage{vmargin} 
\setpapersize{USletter} 
\setmargrb{1in}{1in}{1in}{1in} % --- sets all four margins.  Cool. 
\numberwithin{equation}{section} 
\newtheorem{theorem}{Theorem}[section] 
 
\newtheorem{lemma}[theorem]{Lemma}

\newtheorem{defn}[theorem]{Definition}

\theoremstyle{definition}

\newtheorem{example}[theorem]{Example} 
\newtheorem{remark}[theorem]{Remark}

\newcommand{\R}{\mathbb{R}}
\newcommand{\C}{\mathbb{C}}

\newcommand{\Q}{\mathbb{Q}}
\newcommand{\PP}{\mathbb{P}}

\title{The Convex Hull of a Variety}

\author{Kristian Ranestad}
\address{Kristian Ranestad \\
Matematisk Institutt \\
      Universitetet I Oslo \\
PO Box 1053 \\
Blindern, NO-0316 Oslo, Norway}
\email{ranestad@math.uio.no}
\author{Bernd Sturmfels} 
\address{Bernd Sturmfels\\ 
Department of Mathematics\\ 
      University of California\\ 
      Berkeley, California 94720,
      USA} 
\email{bernd@math.berkeley.edu} 
\begin{document} 

\begin{abstract}
We present a characterization, in terms of projective biduality,
for the hypersurfaces appearing in the boundary of the
convex hull of a compact real algebraic variety.
\end{abstract}

\maketitle

\section{Formula for the Algebraic Boundary}
Convex algebraic geometry is concerned with  the algebraic study of convex sets 
that arise in polynomial optimization.
One topic of recent interest is the convex hull 
${\rm conv}(C)$ of a compact algebraic curve $C$ in $\R^n$.
Various authors have studied 
semidefinite representations \cite{Hen, Sch},  facial structure \cite{SSS, Vin},
and volume estimates \cite{BB, SS} for such convex bodies.
 In \cite{RS} we characterized the boundary geometry 
of ${\rm conv}(C)$ when $n=3$.
The boundary is formed by the edge surface and the tritangent planes,
the degrees of which we computed in \cite[Theorem 2.1]{RS}.
Here, we extend our approach to varieties of any dimension in any $\R^n$.

Throughout this paper, we let $X$ denote a compact algebraic variety
in $\R^n$
which affinely spans $\R^n$. We write $\bar{X}$ for
the Zariski closure of $X$ in complex projective space $\C \PP^n$. 
Later we may add further hypotheses on $X$, e.g., that the complex variety 
$\bar X$ be smooth or irreducible.

The convex hull $P = {\rm conv}(X)$ of $X$ is an $n$-dimensional
 compact convex semialgebraic subset of $\R^n$. We are interested in the
 boundary $\partial P$ of $P$.  Basic results in convexity \cite[Chapter 5]{Gru} and 
 real algebraic geometry \cite[Section 2.8]{BCR} ensure that
  $\partial P$ is a semialgebraic set of pure dimension $n-1$.  The 
  singularity structure of this boundary has been studied by S.D. 
  Sedykh \cite{Sed1,Sed2}.
 Our object of interest is the {\em algebraic boundary}
$\partial_a P$,
by which we mean the Zariski closure of $\partial P$ in $\C \PP^n$.
Thus $\partial_a P$ is a closed subvariety  in $\C \PP^n$  of pure dimension $n-1$.
We represent $\partial_a P$ by the 
 polynomial in $\R[x_1,\ldots,x_n]$ that vanishes on $\partial P$. 
This polynomial is unique up to a multiplicative constant as we require
it to be squarefree. Our ultimate goal is to compute the polynomial
representing the algebraic boundary $\partial_a P $.

We write $X^*$ for the projectively dual variety to $\bar{X}$. The dual variety $X^*$
lives in the dual projective space $(\C \PP^n)^\vee$. It is the Zariski closure
of the set of all hyperplanes that are tangent to $\bar{X}$ at a regular point. According
to the {\em Biduality Theorem} of projective geometry, we have $(X^*)^* = \bar{X}$.
We refer to \cite[\S I.1.3]{GKZ}  for a proof of this important result.

For any positive integer $k $ we let
 $X^{[k]}$ denote the Zariski closure in $(\C \PP^n)^\vee$ of the
set of all hyperplanes that are tangent to $\bar{X}$ at  $k$ regular points
that span a $(k{-}1)$-plane.
Thus $X^{[1]} = X^*$ is the dual variety. We consider the following nested chain
of algebraic varieties:
$$ X^{[n]}\, \subseteq \,\cdots \, \subseteq \,
X^{[2]} \,\subseteq  \, X^{[1]} \,\subseteq \,(\C \PP^n)^\vee. $$
Our objects of interest is the dual variety, back in $\C \PP^n$,
to any $X^{[k]}$ appearing in this chain.

To avoid anomalies, we make the assumption that only finitely many hyperplanes are 
tangent to $\bar{X}$ at infinitely many points.
Without this assumption, the relevant projective geometry
is much more subtle, as seen
in the recent work of Abuaf \cite{Abu}.  
With this assumption,
for small values of $k$, the dual variety  $(X^{[k]})^*$ equals the {\em $k$-th secant variety}
 of $X$, the closure of the union of all $(k-1)$-dimensional linear spaces that intersect $X$ in at least $k$ points.

 The codimension of 
this secant variety is at least $2$ if $\,\,k\leq\lfloor\frac{n}{{\rm dim}(X)+1}\rfloor$.
Let $r(X)$ be the minimal integer $k$ such that the $k$-th secant variety of $X$ has dimension at least $n-1$.  Thus we have $r(X)\geq\lceil\frac{n}{{\rm dim}(X)+1}\rceil$.
The inequality  $k \geq r(X)$ is necessary for  $(X^{[k]})^*$ to be a hypersurface.
The main result in this article is the following formula for the convex hull.

\begin{theorem}   \label{thm:main}
Let $X$ be a smooth and compact real algebraic variety that affinely 
spans $\R^n$, and assume that only finitely many hyperplanes in
$\C\PP^n$ are tangent to the corresponding projective variety $\bar {X}$ at infinitely many points.
The algebraic boundary of its convex hull, $P = {\rm conv}(X)$, 
is computed by biduality as follows:
\begin{equation}
\label{mainformula}
 \partial_a  P \quad \subseteq \quad \bigcup_{k=r(X)}^{n}  (X^{[k]})^*. 
 \end{equation}
In particular, every irreducible component of $ \partial_a P$ is a component of $(X^{[k]})^*$ for some $k$.
\end{theorem}

Since $\partial_a P $ is a hypersurface, at least
one of the $(X^{[k]})^*$ must be a hypersurface.
However, others may have higher codimension and these can be
removed from the union.
The reason for possibly not having equality in (\ref{mainformula}) is that some
hypersurface component in $(X^{[k]})^*$ may have no real points, or its real points may be
disjoint from the boundary of $P = {\rm conv}(X)$. Such components must also be removed
when we compute the algebraic boundary $\partial_a P$.

When the inclusion $X^{[k]} \,\subseteq  \, X^{[k-1]}$ is proper, the former is part of the singular locus of the latter.  
In particular $X^{[k]}$ is in general part of
the $k$-tuple locus of the dual variety $ X^{[1]}=X^*$.
However, the singular locus of $X^*$ will have further components.
  For example, the dual variety of a curve or surface in $\C\PP^3$ 
 has a cuspidal edge defined, respectively,
  by the osculating planes to the curve, and by planes that intersect the surface in a cuspidal curve.

Our presentation is organized as follows.
In Section 2 we discuss a range of examples which illustrate the formula
(\ref{mainformula}). The proof of Theorem \ref{thm:main} is given in Section 3.
We also examine the case when $X$ is not smooth,
and we extend  Theorem \ref{thm:main} to varieties whose real singularities
are isolated.
 Section 4 features additional examples. These highlight the need to
 develop better symbolic and numerical tools for evaluating
the right hand side of (\ref{mainformula}).

\section{First Examples}

\subsection{Polytopes}

Our first example is the case of finite varieties, when ${\rm dim}(X) = 0$.
Here $P = {\rm conv}(X)$ is a full-dimensional 
convex polytope in $\R^n$, and its algebraic boundary
$\partial_a P $ is the Zariski closure of the union of all facets of $P$.
The formula (\ref{mainformula}) specializes to
$$  \partial_a P \quad \subseteq \quad (X^{[n]})^* . $$
Indeed, $X^{[n]} \subset (\C \PP^n)^\vee$ is the finite set of
 hyperplanes that are spanned by $n$
affinely independent points in $X$. Typically, this includes
hyperplanes that do not support $\partial P$, and these should be removed
when passing from $(X^{[n]})^*$ to $\partial_a P$. It is  important to note
that the Zariski closure, used in our definition of the algebraic boundary
$\partial_a P $, depends on the field $K \subseteq \R$ we are working over.
If we take $K = \R$ then $\partial_a P$ is precisely the union of the
facet hyperplanes of $P$. However, if $K$ is the field of definition of $X$,
say $K = \Q$, then $\partial_a P$ usually contains additional hyperplanes
that are Galois conjugate to the facet hyperplanes.

Here is a tiny example that illustrates this arithmetic subtlety.
Let $n=1$ and take $X$ to be the variety of the univariate polynomial
$ x^5 - 3x + 1$. This polynomial is irreducible over $\Q$ and has three real roots.
The smallest root is $\,\alpha = -1.3888...\,$ and the largest root is $\, \beta =  1.2146...$.
Clearly, $P = {\rm conv}(X) $ is the line segment $[\alpha,\beta]$ in $\R^1$.
If we take $K = \R$ then $\,\partial_a P = \{\alpha,\beta\}$, but
if we take  $K = \Q$ then $\partial_a P$ consists of all five complex roots of $f(x)$.

\subsection{Irreducible Curves}
Let $n=2$ and $X$ an irreducible compact curve in $\R^2$ of degree $d \geq 2 $.
Since $X$ is a hypersurface, we have $r(X) = 1$. Suppose that the curve $X$ does not bound
a convex region in $\R^2$.
 The algebraic boundary of the  convex set $ P = {\rm conv}(X)$ consists
of $X$ and the union of all bitangent lines of $X$. In symbols,
$$ \partial_a P \,\,\, \subseteq \,\,\,  (X^{[1]})^* \,\cup \, (X^{[2]})^* \,\, = \,\, X \,\, \cup \,\, (X^{[2]})^*. $$
For a smooth curve $X$ of degree $d$, the classical Pl\"ucker formulas imply that 
the number of (complex) bitangent lines equals
$(d-3)(d-2)d(d+3)/2$.  Hence, $\partial_a P$ is a curve of degree 
$$ {\rm deg}(\partial_a P)\, \, \,\leq\,\,\,
d \,+ \, \frac{(d-3)(d-2)\,d\,(d+3)}{2}. $$

Next consider the case
where $n = 3$, ${\rm dim}(X)= 1$, and $r(X)=2$. If $X$ is irreducible then  
$$  \partial_a P \,\,\, \subseteq \,\,\,  (X^{[2]})^* \,\cup \, (X^{[3]})^*. $$
The first piece $(X^{[2]})^*$ is the {\em edge surface} of $X$,
and the second piece $(X^{[3]})^*$ is the union of all
{\em tritangent planes}. For a detailed study of this situation, including
pretty pictures of $P$,  and
a derivation of degree formulas for  $(X^{[2]})^*$ and $(X^{[3]})^*$,
we refer to our earlier paper \cite{RS}.
Further examples of space curves are found in
Subsection 4.1 below and in \cite[Section 4]{RS}.

Sedykh and Shapiro \cite{SS} studied {\em convex curves} $X \subset \R^n$
where $n=2r$ is even. Such a curve has the property that 
$|X \cap H| \leq n $ for every real hyperplane $H$. The algebraic boundary 
of a convex curve is the hypersurface of all secant $(r-1)$-planes. In symbols,
$\,\partial_a P \, = \, (X^{[r]})^*$.

\subsection{Surfaces in 3-Space}

Let $X$ be a general smooth compact surface of degree $d$ in $\R^3$.
Confirming classical derivations by Cayley, Salmon and Zeuthen 
\cite[p.313-320]{Sal}, work on enumerative
geometry in the 1970s by
Piene \cite[p.231]{Pie} and Vainsencher \cite[p.414]{Vai} establishes the following
formulas for the degree of the curve $X^{[2]}$, its dual surface
$(X^{[2]})^*$, and the finite set
$ X^{[3]}$ in $(\C \PP^3)^\vee$:
\begin{eqnarray*}
{\rm deg} (X^{[2]}) & = & \frac{d(d-1)(d-2)(d^3-d^2+d-12)}{2} ,\\
{\rm deg}\bigl((X^{[2]})^*\bigr) & = & d(d-2)(d-3)(d^2+2 d-4), \\
& & \\
{\rm deg} (X^{[3]}) 
& = & {\rm deg}\bigl((X^{[3]})^*\bigr) \\
&  =& \frac{d^9-6d^8+15d^7-59d^6
+204 d^5-339 d^4+770 d^3-2056d^2 +1920d}{6}.
\end{eqnarray*}
We can expect the degree of $\partial_a P$ to be bounded above by
$d$ plus the sum of the last two expressions, since
$$ \partial_a P \,\, \subseteq \,\, 
    (X^{[1]})^* \, \cup \, (X^{[2]})^* \, \cup \, (X^{[3]})^* \,\, = \,\,
 X \, \cup \, (X^{[2]})^* \, \cup \, (X^{[3]})^* , $$
 unless $X$ is convex or otherwise special.
For a numerical example consider the case $d=4$, where
we take $X$ to be a compact but non-convex smooth quartic surface in $\R^3$.
The above formulas reveal that the degree of the algebraic boundary
$\partial_a P$ can be as large as
$$ {\rm deg}(X) \,+ \,
{\rm deg}((X^{[2]})^*) \, + \,
{\rm deg}((X^{[3]})^*) \,\, = \,\,
 4 + 160 + 3200 \,\, = \,\, 3364.
$$

\subsection{Barvinok-Novik curve}

We examine the first non-trivial instance of the family of {\em Barvinok-Novik curves}
studied in \cite{BN, Vin}. This is the curve $X \subset \R^4$ parametrically given by
$$ (c_1,c_3,s_1,s_3) \, = \, \bigl(
{\rm cos}(\theta), {\rm cos}(3 \theta),
{\rm sin}(\theta), {\rm sin}(3 \theta) \bigr). $$
We change to complex coordinates by setting $x_j = c_j + \sqrt{-1} \cdot s_j$ and
${\bar x}_j = c_j - \sqrt{-1} \cdot s_j$.
The convex body $P = {\rm conv}(X)$ is the projection of the
$6$-dimensional {\em Hermitian spectrahedron}
$$ \biggl\{ (c_1,c_2,c_3,s_1,s_2,s_3) \in \R^6 \,:\,
\begin{pmatrix} 1 & x_1 & x_2 & x_3 \\ {\bar x}_1 & 1 & x_1  & x_2 \\
{\bar x}_2 & {\bar x}_1 & 1  & x_1 \\ {\bar x}_3 & {\bar x}_2 & {\bar x}_1 & 1 \end{pmatrix} \,
\hbox{is positive semidefinite} \biggr\} $$
under the linear map
$\,\R^6 \rightarrow \R^4, \,(c_1,c_2,c_3,s_1,s_2,s_3) \mapsto (c_1,c_3,s_1,s_3)$.
The curve $X$ is the projection of the curve in $\R^6$ 
that consists of the above Toeplitz matrices that have rank~$1$.

The convex body $P = {\rm conv}(X)$ was studied 
 in \cite[Example 5.5]{SSS}. It is the $4$-dimensional representative of
 the   {\em Barvinok-Novik orbitopes} (cf.~\cite{BN, Vin}).
  Its algebraic boundary equals
 $$ \partial_a  P \,= \, (X^{[2]})^* \,\cup \,(X^{[3]})^*  . $$
The threefold $(X^{[2]})^*$ represents the
$2$-dimensional family of edges of $P$, while
the threefold $(X^{[3]})^*$ represents the
$1$-dimensional family of triangles in $\partial P$,
both of which are described in \cite[Thm.~4.1]{BN}; see also \cite{Vin}.
The defining polynomials of these two hypersurfaces in $\R^4$ are
$$ \!\!\!\!\!\! \begin{matrix} 
 & \bigl\langle\,
     x_3^2  {\bar x}_1^6-2 x_1^3 x_3  {\bar x}_1^3  {\bar x}_3+x_1^6  {\bar x}_3^2+4 x_1^3  {\bar x}_1^3
   -6 x_1 x_3  {\bar x}_1^4-6 x_1^4   {\bar x}_1  {\bar x}_3+12 x_1^2 x_3  {\bar x}_1^2  {\bar x}_3 \qquad&
  \\   \qquad & - \, 2 x_3^2  {\bar x}_1^3  {\bar x}_3    -2 x_1^3 x_3  {\bar x}_3^2                   
   -3 x_1^2  {\bar x}_1^2+4 x_3  {\bar x}_1^3+4 x_1^3  {\bar x}_3-6 x_1 x_3  {\bar x}_1  {\bar x}_3
   +x_3^2  {\bar x}_3^2 \,\bigr\rangle&
\end{matrix}  \hbox{from} \,\, (X^{[2]})^* , $$
$$ \qquad \qquad  \qquad \quad
\hbox{and} \qquad \qquad \qquad \bigl\langle \, x_3 {\bar x}_3 \,-\, 1 \,\bigr\rangle 
\,\,\, = \,\,\,\bigl\langle c_3^2 \, + \,s_3^2 \,-\, 1 \bigr\rangle
 \qquad \qquad \qquad \quad \, 
\hbox{from} \,\,\, (X^{[3]})^* .
\qquad \qquad  $$
Both of these threefolds are irreducible components of the 
ramification locus that arises
when we project the hypersurface of singular Toeplitz matrices from $\R^6$ into $\R^4$
as above.

\section{Proof of the Formula}

 We turn to the proof of our biduality formula for the algebraic boundary of $P = {\rm conv}(X)$.
 
 \begin{proof}[Proof of Theorem \ref{thm:main}]
We first prove that the supporting hyperplane of any exposed face $F$ of $P$ lies in $X^{[k]}$  for some $k$. 
Suppose that ${\rm dim}(F) = k-1$ and let  $L_F$ be the projective span of $F$.
By Carath\'eodory's Theorem, every point of $F$ lies in the convex hull of $k$ distinct points 
on $X$.  In particular, the $(k-1)$-plane $L_{F}$ intersects $X$ in at least $k$ points 
that span a $(k-1)$-simplex in $F$.
 If $H$ is a supporting hyperplane for $F$, then $H$ contains $F$ and is the boundary of a halfspace that contains $X$.   
 Since $X$ is smooth, the tangent plane to $X$ at each point $q\in 
 X\cap F\subseteq X\cap H$ must therefore be contained in $H$.
We conclude that  $\,[H]\in X^{[k]}$.

Now, consider any irreducible hypersurface $Y \subset \C \PP^n$ whose real locus has
full-dimensional intersection with the boundary $\partial P \subset \R^n$.
We need to show that
$Y$  is a component of $(X^{[k]})^*$ for some $k$.   
In the next paragraph we give an overview of the proof that follows thereafter.

First, we shall identify the relevant number $k=k_{Y}+1$, by 
the property that $Y$ has a linear space 
of dimension $k_{Y}$ through every point.
In fact, we shall prove that $Y$ contains a unique $k_{Y}$-plane through a general point of $Y$.
Thus, at a general point, the hypersurface $Y$ is locally a fibration. 
In particular, the general point in $\partial P\cap Y$ lies in a 
$k_{Y}$-plane that intersects $P$ along a $k_{Y}$-dimensional face.
Subsequently, we will show that the supporting hyperplanes of these faces 
are tangent to $Y$ along these $k_{Y}$-planes, before we prove that $Y^*\subseteq X^{[k]}$.
From this, we shall finally conclude that $Y$ is a component of $(X^{[k]})^*$.

Let $q$ be a general smooth point in the $(n-1)$-dimensional semialgebraic set
$\partial P\cap Y$. Since the union of the exposed faces of $P$ is dense in $\partial P$,
there exists an exposed face $F_q$  that has $q$ in its
relative interior.    
 The hypersurface $Y$ contains the boundary of $P$ locally at $q$, and hence it contains the face $F_q$. 
 Since $Y$ is a variety, it contains the projective span $L_{F_q}$ of the face $F_q$.  Let $k_{Y}={\rm dim}(L_{F_q})$. Since
 $q$ is a general smooth point in $\partial P\cap Y$,
the hypersurface $Y$ contains a $k_{Y}$-plane through every point of~$Y$. 
 In fact, since $F_{q}$ is an exposed face, it is the unique face 
 through $q$, so $Y$ contains a unique $k_{Y}$-plane through every 
 % Kristian, please check (Bernd 04/12)
 general
 point. 

Next, let $H$ be a hyperplane that exposes the $k_{Y}$-dimensional 
face $F_{q}$ of $P$.  We will show that $H$ coincides with the tangent hyperplane $H_q$ to $Y$ at $q$.
  As $q$ is a general interior point in $F_{q}$, we then conclude 
  that $H$ is tangent to $Y$  along the entire $k_{Y}$-plane $L_{F_{q}}$. 
   The key to our argument is that $H$
   is assumed to be tangent to $X$ at the points $X\cap F_{q}$ that span $L_{F_{q}}$.

%Fix a hyperplane $H$ as above.
If $Y=L_{F_{q}}$ is itself a hyperplane, there is nothing to prove, except to note that $k_{Y}=n-1$, that  $H=H_{q}$, and that $Y^*$ is an isolated point in $X^{[n]}$.
Otherwise, we compare $H$ and the tangent plane $H_{q}$ via a local parameterization of $Y$ at $q$.  Let $k=k_{Y}+1$ and $m={\rm dim}(X)$, let $p_{1},...,p_{k}$ be points in $X\cap F_{q}$ that 
affinely span $F_{q}$, and let 
$$\qquad \gamma_{i}: t_{i}=(t_{i,1},...,t_{i,m})\mapsto(\gamma_{i,1}(t_{i}),...,\gamma_{i,n}(t_{i}))
\qquad \qquad \hbox{(for $i=1,...,k$)} $$
be local parameterizations of $\bar{X}$ at the points $p_{i}$.
% in $\R^n$.
The point $q$ lies in the affine-linear span of the points $p_{i}$, so 
$q=\sum_{i}^k a_{i}p_{i}$ for some real coefficients $a_{i}$ with $\sum a_i = 1$.
There may be polynomial relations in the local parameters $t_{i}$ 
defining $k$-tuples of points 
in $X$ whose affine-linear span lies in $Y$. These relations define 
a subvariety $Z$ in the Cartesian product $\bar{X}^k$ that contains the 
$k$-tuple $(p_{1},\ldots,p_{k})$.  
A local parametrization of $L_{F_{q}}$ at $q$ has the form 
$$\alpha:  u=(u_{1},...,u_{k_{Y}}) \mapsto (\alpha_{1}(u), ..., \alpha_{n}(u))$$ 
with affine-linear functions $\alpha_{i}$ in the $u_{i}$.
Since $Y$ is locally a fibration, the 
algebraic functions $\gamma_{i}$ and $\alpha$ provide a local 
parameterization of the complex variety $Y$ at the point $q$:
$$\begin{matrix}
\Gamma:& \C^{k_{Y}}\times Z&\qquad\to&\qquad \C^n\qquad\qquad\qquad\qquad\qquad\qquad\qquad\\
&(u, t_{1},...,t_{k})&\qquad\mapsto&\qquad \alpha(u)+\sum_{i=1}^k 
a_{i}(\gamma_{{i}}(t_{i}))+\epsilon(u,t_{1},\ldots,t_{k})\\
\end{matrix}
$$
Here, the function
 $\epsilon$ only contains terms of order at least two in the parameters.
 The tangent space $H_{q}$ at $q$ is spanned by the linear terms in the above parameterization.  But these linear terms lie in the span 
of $(\alpha_{1}(u),...,\alpha_{n}(u))$ and the linear terms in 
$(\gamma_{1},...,\gamma_{k})$.  The former span $L_{F_{q}}$, 
while the latter span the tangent spaces to $\bar{X}$ at each of the points 
$p_{i}$.  So, by assumption they all lie in the hyperplane $H$ that 
supports $\partial P$ at  $F_q$.
 Therefore, the hyperplane $H$ must coincide with the tangent plane $H_{q}$ to $Y$ at $q$.
 Since $q$ is a general point not just in $Y$ but also in $L_{F_{q}}$, 
 we conclude that $H$ is tangent to $Y$ along the entire plane $L_{F_{q}}$.

 We have shown that the tangent hyperplanes to $Y$ are constant along the $k_{Y}$-planes contained in $Y$.
 This implies that  the dimension of the dual variety $Y^*$ is equal to $n-k$
 where   $k=k_{Y}+1$. Locally around the point $q$, these
 tangent hyperplanes support faces of dimension $k_{Y}=k-1$
  the convex body $P$. This ensures that the inclusion
$Y^*\subseteq X^{[k]}$ holds.
  
We next claim that $Y^*$ is in fact an irreducible component
of the variety $X^{[k]}$. This will be a consequence of the following 
general fact which we record as a lemma.

\begin{lemma} Every irreducible component of  $X^{[k]}$ has dimension at most $n-k$.  
\end{lemma}
\begin{proof} 
Let  $W\subseteq X^{[k]}$ be a component, and let 
    $k_{W}$ be the minimal $l$ such that $W$ is not contained in 
    $X^{[l+1]}$.  Then $k_{W}\geq k$ and $W$ is a component of 
    $X^{[k_{W}]}$. 
   
    Let $CX\subset \C \PP^n\times (\C \PP^n)^\vee$ be the conormal variety of $\bar{X}$, 
    the closure of the set of pairs $(p,[H])\in \C \PP^n\times (\C \PP^n)^\vee$ 
    such that the hyperplane $H$ is tangent at the smooth point $p\in X$.
     It has dimension $n-1$.  By assumption, the projection $\rho:CX\to (\C \PP^n)^\vee $
      into the dual space has only finitely many infinite fibers.
    Therefore $X^*=\rho(CX)$ is a hypersurface and  $W$ is part of its $k_{W}$-tuple locus.  
    If $[H]$ is a general point in $W$, then $X^*$ has at least $k_{W}$ branches at $[H]$. 
    Let $(p_{1},[H]),...,(p_{k_{W}},[H])$ be smooth points in $CX$ in the fiber over $[H]$, such that $p_{1},...,p_{k_{W}}$ are linear independent points on $X$.  
    Consider the tangent spaces $T_{1},...,T_{k_{W}}$ to $CX$ at these points, and let $U_{i}=\rho_{T}(T_{i}), i=1,...,k_{W}$ be the corresponding linear spaces 
    in the tangent cone to $X^*$ at $[H]$, where $\rho_{T}$ is the map induced by $\rho$ on tangent spaces.
    Then the intersection $U_{1}\cap \cdots\cap U_{k_{W}}$ contains the tangent space to $W$ at $[H]$.  But  $p_i \in U_{i}^{\bot}$, so the orthogonal complement of
    the intersection satisfies
           $$(U_{1}\cap \cdots \cap U_{k_{W}})^{\bot}
           \,\,\, = \,\,\,
           {\rm span}(U_{1}^{\bot}\cup ...\cup U_{k_{W}}^{\bot})
           \,\,\, \supseteq \,\,\, {\rm span}( p_{1},...,p_{k_{W}}).$$
           We conclude that the plane
 $U_{1}\cap \cdots \cap U_{k_{W}}$ has codimension at least $k_{W}$ at $[H]$,
 and therefore the variety
 $W$ has codimension at least $k_W$ in $(\C \PP^n)^\vee$.
   Since $k_{W}\geq k$ the lemma follows.
\end{proof} 

At this point, we are pretty much done.
To recap, recall that we have shown $Y^* \subseteq X^{[k]}$,
${\rm dim}(Y^*) =n-k$ and ${\rm dim}(X^{[k]}) \leq n-k$.
 If $X^{[k]}$ is irreducible, then we have
$Y^* = X^{[k]}$ and $Y = (X^{[k]})^*$ follows. Otherwise, if $X^{[k]}$ has several components,
then its dual $(X^{[k]})^*$ is the union of the dual varieties of each component.  
One of these components is $Y$, and hence $Y^*$ is a component of $(X^{[k]})^*$.
Therefore, the formula (\ref{mainformula}) 
in Theorem \ref{thm:main} is indeed true.
\end{proof}

\medskip
Theorem \ref{thm:main} extends in a straightforward manner
to reduced and reducible compact real algebraic sets with isolated singularities.
A colorful picture of a trigonometric space curve $X$  with  a singularity on the boundary of
$P = {\rm conv}(X)$ is shown in \cite[Figure 6]{RS2}. Also,
in Subsection 4.1 below we shall examine a reducible space curve with isolated
singularities with the property that some (finitely many) hyperplanes that are tangent
at infinitely many points.

Let $X\subset\R^n$ be a finite union of compact varieties, and assume that $X$ has only isolated singularities.  As before, we write $\bar{X}$ be its Zariski closure in $\C\PP^n$.
For any positive integer $k $ we now take
 $X^{[k]}$ to be the Zariski closure in $(\C \PP^n)^\vee$ of the
set of all hyperplanes that are tangent to $\bar{X}$ at  $k-s$ regular points and pass through $s$ singularities on $X$, for some $s$, such that the $(k-s)+s=k$ points span a $(k{-}1)$-plane.
Thus $X^{[1]}$ contains the dual variety, but, in addition, it 
also contains a hyperplane for each isolated singularity of $X$. We consider, as above, 
the nested chain of projective varieties
$$ X^{[n]}\, \subseteq \,\cdots \, \subseteq \,
X^{[2]} \,\subseteq  \, X^{[1]} \,\subseteq \,(\C \PP^n)^\vee. $$
The algebraic boundary of $P={\rm conv}(X)$ is dual to the various
 $X^{[k]}$ appearing in this chain:

\begin{theorem}   \label{thm:sing}
Let $X$ be a finite union of compact real algebraic varieties that affinely spans $\R^n$, 
and assume that $X$ has only isolated singularities and that only finitely many hyperplanes
in $\C\PP^n$ are 
tangent to $\bar {X}$ at infinitely many points.
The algebraic boundary of its convex hull, $P = {\rm conv}(X)$, 
is computed by biduality using the same formula (\ref{mainformula})
as in Theorem \ref{thm:main}.
In particular, every irreducible component of $ \partial_a P$ 
is a component of $(X^{[k]})^*$ for some~$k$.
\end{theorem}

\begin{proof}  Following the argument of the proof of 
Theorem \ref{thm:main}, we first note that a hyperplane $H$ that supports 
a $(k-1)$-dimensional face of $P$ must intersect $X$ in $k$ points that span the face.  Furthermore, $H$ must be tangent to $X$ at the smooth
intersection points. Let $Y$ be an irreducible component having full-dimensional 
intersection with the boundary $\partial P$ of ${\rm conv}(X)$, and let $q$ be a
general smooth point on $\partial P\cap Y$.
In the notation of the above proof, a local parameterization of $Y$ at $q$ will 
involve singular points $p_{1},...,p_{s}$ and smooth points 
$p_{s+1},...,p_{k}$.  The $k$-tuples $(p_{1},\ldots,p_{k})$ of  
points whose linear span is contained in $Y$ form a subvariety $Z$ in 
the Cartesian product $X^k$.  Since the singular points are isolated, 
we may assume that the restriction of $Z$ to the first $s$ 
factors is a point. The hypersurface $Y$ is therefore a cone with vertex containing 
the $s$ singular points. The tangent hyperplane to $Y$ at $q$ contains the vertex 
and the tangent spaces at the  $k-s$ smooth points, 
so it coincides with the supporting hyperplane $H$.
The latter part of the proof of  Theorem \ref{thm:main} applies directly to arrive
at the same conclusion.
\end{proof}

At present, we do not know how to extend our formula 
(\ref{mainformula}) for the algebraic boundary to the
convex hull of a compact real variety $X$ whose real singular locus has dimension
$\geq 1$. Also, we do not yet know how to remove the hypothesis
that only finitely many hyperplanes are 
tangent to $\bar{X}$ at infinitely many points.
This issue is related to the study of degeneracies in \cite{Abu}
and we hope that the techniques introduced in that paper
will help for our problem.

\section{More Examples and Computational Thoughts}

We further illustrate our formula for the algebraic boundary of the convex hull
of a real variety with three concrete examples, starting
with a curve that is reducible and singular.

\subsection{Circles and spheres in 3-Space}

 Let $ n = 3$ and suppose that $X = C_1 \cup C_2 \cup 
\cdots \cup C_r $ is the reducible (and possibly singular)
curve obtained by taking the union of
 a collection of $r \geq 3$ sufficiently general circles $\,C_i\,$ that lie in 
various planes in $\R^3$. We have
\begin{equation}
\label{CircIn3}
 \partial_a P \,\, \subseteq \,\,  (X^{[2]})^* \, \cup \, (X^{[3]})^* .
 \end{equation}
The surface $(X^{[3]})^*$ 
is the union of planes that are tangent to three of the circles and planes spanned by the circles.
The edge surface $(X^{[2]})^*$ decomposes into quadratic surfaces, namely,
its components are cylinders formed by
stationary bisecant lines defined by pairs of circles.

\begin{figure}
  \centering
  \includegraphics[width=0.59\textwidth]{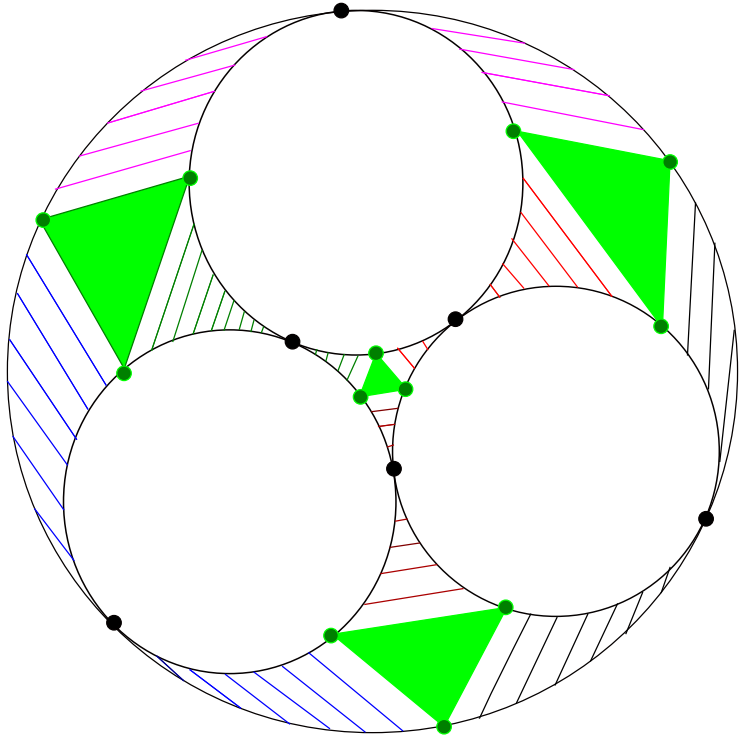}
  \caption{Schlegel diagram of the convex hull of four pairwise touching circles}
  \label{fig:Schlegel4circles}
\end{figure}

For a concrete configuration, consider the convex hull of $r= 4 $ 
pairwise touching circles in $3$-space.   
The surface $ (X^{[2]})^*$ is a union of six cylinders, each wrapped around two of the circles, while $(X^{[3]})^*$ is the union of
 planes tangent to three of the circles (four of which contain the fourth circle). A picture of this
$3$-dimensional convex body $P$ is shown in
Figure~\ref{fig:Schlegel4circles}.
 Its boundary consists of $6+(4+4)= 14$ distinct surface
patches, corresponding to the pieces in (\ref{CircIn3}), which holds with equality.
There are six cylinders, four planes touching exactly three of the circles,
 and four planes containing one of the circles and touching the three others.

A nice symmetric
representation of the curve $X = C_1 \cup C_2 \cup C_3 \cup C_4$ is given by the ideal
$$ \langle \,a\, c \, g \, t \,\,,\,\, a^2 + c^2 + g^2 + t^2 - 2 ac - 2ag - 2at - 2c
g - 2 ct - 2gt \,\rangle, $$
where the variety of that ideal is to be taken  inside the probability simplex
$$\Delta_3 \,\, = \,\, \{ \,(a,c,g,t) \in \R^4_{\geq 0} \,: \,a+c+g+t = 1 \,\}. $$
The convex body $P$  looks combinatorially like a $3$-polytope with $18$ vertices, $36$ edges and $20$ cells. Eight of the
$20$ cells are  flat facets. First, there are the planes of the
circles themselves. For instance, the facet in the plane $t = 0$ is the
disk $\{a^2 + c^2 + g^2 \leq 2 ac + 2 ag + 2cg\}$ in the triangle $ \{\,a+c+g = 1 \,\}$.
 Second, there are
four triangle facets, formed by the unique planes that are tangent to 
exactly three of the  circles. The equations of these facet planes are
$$
\begin{matrix}
     P_a \,=\,   -a+2 c+2 g+2 t  , & &&
     P_c \,=\,   2 a-c+2 g+2 t  , \\
     P_g \,\,=\,\,   2 a+2 c-g+2 t  , & & &
     P_t  \,=\,   2 a+2 c+2 g-t  .
  \end{matrix}
$$
The remaining $12$ cells in $\partial P$ are 
quadratic surface patches that arise from the pairwise convex hull
of any two of the four circles. This results in $6$ quadratic surfaces
each of which contributes two triangular cells to the boundary. The
equations of these six surfaces are
$$ 
\begin{matrix}
Q_{ac} \,=\,   \,a^2 + c^2 + g^2 + t^2 + 2( ac - ag - cg - at - ct - gt ) \,  , \\
Q_{ag} \,=\,   \,a^2 + c^2 + g^2 + t^2 - 2(ac - ag + cg + at + ct + gt ) \,  , \\
Q_{ag} \,=\,  \, a^2 + c^2 + g^2 + t^2 - 2(ac + ag - cg + at + ct + gt ) \,  , \\
Q_{cg} \,=\,   \, a^2 + c^2 + g^2 + t^2 - 2(ac + ag + cg - at + ct + gt ) \,  , \\
Q_{ct} \,=\,   \, a^2 + c^2 + g^2 + t^2 - 2(ac + ag + cg + at - ct + gt ) \,  , \\
Q_{gt} \,=\,   \, a^2 + c^2 + g^2 + t^2 - 2(ac + ag + cg + at + ct - gt ) \,  .\\
\end{matrix}
$$

Each circle is subdivided into six arcs of equal length. Three of the nodes arise
from intersections with other circles, and the others
are the intersections with the planes $P_a, P_c, P_g, P_t$.
This accounts for all $18$ vertices and $24$ ``edges'' that are  arcs.
The other $12$ edges of $\partial P$ are true edges: they arise from
the four triangles. These are shown in green in the \emph{Schlegel
diagram} of Figure~\ref{fig:Schlegel4circles}.
The 12 cells
corresponding to the six quadratic surfaces are the 12 ruled
cells in the diagram, and they come in pairs according to the six
different colors. The six intersection points among the 4 circles are
indicated by black dots, whereas the remaining twelve vertices
correspond to the green dots which are vertices of our four green triangles.

\begin{figure} \centering
  \includegraphics[width=0.5\textwidth]{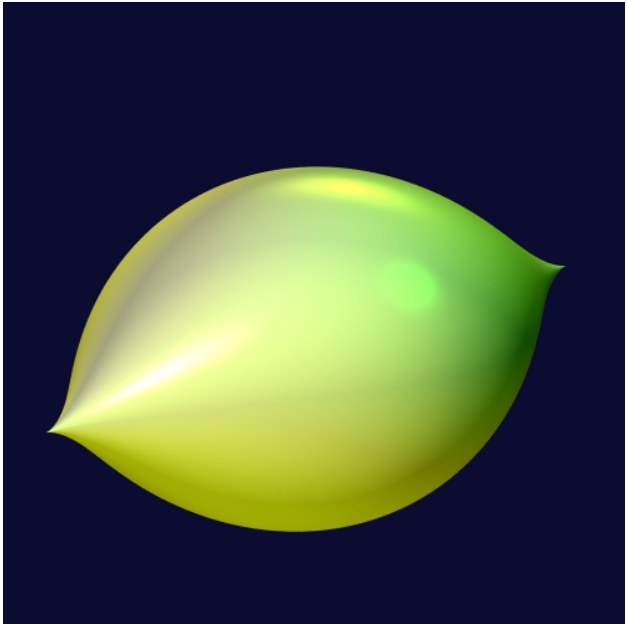}
  \caption{The Zitrus surface $\, x^2+z^2+(y^2-1)^3 = 0 $}
    \label{fig:zitrus}
\end{figure}

\subsection{Zitrus}

We have seen that the convex hull of algebraic surfaces in $\R^3$
can have surfaces of very high degree in its boundary. For instance,
if $X$ is a general smooth surface of degree $d=6$ then
the bitangent surface  $ \,(X^{[2]})^* \,$ has degree $\, 3168$. On
the other hand, that number can be expected to drop substantially
for most singular surfaces. Let us consider the sextic
$$ f(x,y,z) \quad = \quad x^2+z^2+(y^2-1)^3 . $$
The surface $X = V(f)$ in $\R^3$ is taken from Herwig Hauser's beautiful
{\em Gallery of Algebraic Surfaces}.
The name given to that surface is {\em Zitrus}.
It appears on page 42-43 of the catalog \cite{imag} of the exhibition {\em Imaginary}.
For an electronic version see  {\tt www.freigeist.cc/gallery.html}.

We choose affine coordinates $(a,b,c)$ on the space of planes
$a x + b y + c z + 1 = 0$ in $\R^3$. In these coordinates, the
variety $X^{[2]}$ is the union of two quadratic curves given by the ideal
$$ \langle b+1,27 a^2+27 c^2-16 \rangle \,\cap \, 
\langle b-1,27 a^2+27 c^2-16 \rangle. $$
These curves parametrize the tangent planes that pass through one of the
two singular points of the Zitrus. Each curve dualizes to a singular quadratic surface,
and  $(X^{[2]})^*$ is given by
$$ \langle 16x^2-27y^2+16z^2+54y-27 \rangle \,\cup \,
\langle 16x^2-27y^2+16z^2-54y-27 \rangle.
$$
The Zitrus $X$ has no tritangent planes, so 
$\partial_a P =  X \,\cup \, (X^{[2]})^*$, and we conclude that the
algebraic boundary of the
{\em  convexified Zitrus}  $\,P  = {\rm conv}(X)$  has degree $10 = 6+2+2$.

We now perturb the polynomial $f$ and  consider the smooth surface
$\tilde X = V(\tilde f)$ defined by
$$ \tilde f(x,y,z)  \quad = \quad x^2+z^2+(y^2-1)^3  - 1. $$
The curve of bitangent planes, ${\tilde X}^{[2]}$, has again two components.
It is defined by the ideal
$$ \begin{matrix} \langle\, b \,,\,a^2+c^2-1 \,\rangle \, \cap \, \langle \,
90 a^2 b^2-96 b^4+90 b^2 c^2-129 a^2+128 b^2-129 c^2+48 \,, \qquad \\ 
\qquad \qquad \qquad \qquad \qquad 135 a^4-144 
       b^4+270 a^2 c^2+135 c^4-6 a^2+272 b^2-6 c^2-48 \,\rangle.
       \end{matrix}  $$
The first curve dualizes to the cylinder $\{x^2+z^2=1\}$. The other component of
the boundary surface $({\tilde X}^{[2]})^*$ has degree $16$. Its defining polynomial has
$165$ terms which start as follows:
$$ 16777216 x^{16}-169869312 x^{14}y^2+1601372160 x^{12} y^4-7081205760
x^{10}   y^6+26435102976 x^8 y^8 - \cdots $$

\subsection{Grassmannian}
\label{subs:grass}

We consider the oriented Grassmannian $X = {\rm Gr}(2,5)$ of oriented
two-dimensional linear subspaces of $\R^5$. This is the $6$-dimensional
subvariety of $\R^{10}$ defined~by
$$
\begin{matrix}
\langle \,
p_{12}^2 +
p_{13}^2 +
p_{14}^2 +
p_{15}^2 +
p_{23}^2 +
p_{24}^2 +
p_{25}^2 +
p_{34}^2 +
p_{35}^2 +
p_{45}^2 -1, \,
p_{12} p_{34} - p_{13} p_{24} + p_{14} p_{23} , \quad \\
p_{12} p_{35} {-} p_{13} p_{25} {+} p_{15} p_{23} , \,
p_{12} p_{45} {-} p_{14} p_{25} {+} p_{15} p_{24} , \,
p_{13} p_{45} {-} p_{14} p_{35} {+} p_{15} p_{34} , \,
p_{23} p_{45} {-} p_{24} p_{35} {+} p_{25} p_{34}  \,
\rangle .
\end{matrix}
$$
Its convex hull $P = {\rm conv}(X)$ is a {\em Grassmann orbitope},
a class of convex bodies that are of interest to differential geometers.
We refer to \cite{Mor}, \cite[\S 7]{SSS}, and the references given therein.
The determinant of the Hermitian matrix in the spectrahedral representation of
$P$ in \cite[Theorem 7.3]{SSS} has degree $8$ and it factors into
two quartic factors. Only one of these two factors is relevant for us,
and we display it below.
Namely, the algebraic boundary
$\,\partial_a P = (X^{[4]})^*\,$ is the irreducible hypersurface of  
degree $4$ represented by the polynomial
$$
\begin{matrix}
  p_{12}^4 +p_{13}^4 +p_{14}^4 +p_{15}^4 +p_{23}^4
+p_{24}^4 +p_{25}^4 +p_{34}^4 +p_{35}^4 +p_{45}^4 \\
+2 p_{12}^2 p_{13}^2 +2 p_{12}^2 p_{14}^2 +2 p_{13}^2 p_{14}^2
+2 p_{12}^2 p_{15}^2
+2 p_{13}^2 p_{15}^2 +2 p_{14}^2 p_{15}^2
+2 p_{12}^2 p_{23}^2 +2 p_{13}^2 p_{23}^2 -2 p_{14}^2 p_{23}^2 \\
-2 p_{15}^2 p_{23}^2 +2 p_{12}^2 p_{24}^2 -2 p_{13}^2 p_{24}^2
+2 p_{14}^2 p_{24}^2 -2 p_{15}^2 p_{24}^2
+2 p_{23}^2 p_{24}^2
+2 p_{12}^2 p_{25}^2 -2 p_{13}^2 p_{25}^2 -2 p_{14}^2 p_{25}^2 \\
+2 p_{15}^2 p_{25}^2 +2 p_{23}^2 p_{25}^2 +2 p_{24}^2 p_{25}^2
-2 p_{12}^2 p_{34}^2 +2 p_{13}^2 p_{34}^2 +2 p_{14}^2 p_{34}^2
-2 p_{15}^2 p_{34}^2 +2 p_{23}^2 p_{34}^2 +2 p_{24}^2 p_{34}^2 \\
-2 p_{25}^2 p_{34}^2 -2 p_{12}^2 p_{35}^2 +2 p_{13}^2 p_{35}^2
-2 p_{14}^2 p_{35}^2 +2 p_{15}^2 p_{35}^2 +2 p_{23}^2 p_{35}^2
-2 p_{24}^2 p_{35}^2
+2 p_{25}^2 p_{35}^2 +2 p_{34}^2 p_{35}^2  \\
-2 p_{12}^2 p_{45}^2 -2 p_{13}^2 p_{45}^2 +2 p_{14}^2 p_{45}^2
+2 p_{15}^2 p_{45}^2 -2 p_{23}^2 p_{45}^2 +2 p_{24}^2 p_{45}^2
+2 p_{25}^2 p_{45}^2
+2 p_{34}^2 p_{45}^2
  +2 p_{35}^2 p_{45}^2 \\
+8 p_{13} p_{14} p_{23} p_{24}
+8 p_{13} p_{15} p_{23} p_{25}
+8 p_{14} p_{15} p_{24} p_{25}
-8 p_{12} p_{14} p_{23} p_{34}
+8 p_{12} p_{13} p_{24} p_{34} \\
-8 p_{12} p_{15} p_{23} p_{35}
+8 p_{12} p_{13} p_{25} p_{35}
+8 p_{14} p_{15} p_{34} p_{35}
+8 p_{24} p_{25} p_{34} p_{35}
-8 p_{12} p_{15} p_{24} p_{45} \\
+8 p_{12} p_{14} p_{25} p_{45}
-8 p_{13} p_{15} p_{34} p_{45}
-8 p_{23} p_{25} p_{34} p_{45}
+8 p_{13} p_{14} p_{35} p_{45}
+8 p_{23} p_{24} p_{35} p_{45} \\
-2 p_{12}^2-2 p_{13}^2-2 p_{14}^2-2 p_{15}^2-2 p_{23}^2
-2 p_{24}^2-2 p_{25}^2-2 p_{34}^2-2 p_{35}^2-2 p_{45}^2
+1.
\end{matrix}
$$
This hypersurface represents a $6$-dimensional family of $3$-dimensional facets of $P$.
Each facet of $P$ is a  $3$-dimensional ball. It meets
the variety $X$ in its boundary, which is a $2$-sphere.

\subsection{Computing}

This paper raises the following algorithmic problem:
given a projective variety $X$, either by its ideal or by a parametrization,
how to compute the equations defining $(X^{[k]})^*$ in practise?
The passage from  $X$ to $X^{[k]}$ can be phrased as an elimination problem
in a fairly straightforward manner. In principle, we can use any
Gr\"obner-based computer algebra system to perform that elimination task.
However, in our experience, this approach only succeeds for 
tiny low-degree instances. Examples such as the Grassmannian in
Subsection \ref{subs:grass} appear to be out of reach for a general purpose 
implementations of our formula (\ref{mainformula}).

 Even the first instance $k=1$, which is the passage 
from a variety $X$ to its dual variety $X^*$, poses a considerable
challenge for  current computational algebraic geometry software.
The case of plane curve is still relatively easy, and it has been addressed
in the literature \cite{BE}.
However, what we need here is the case when $X$ is not a hypersurface
but $X^*$ is. The first interesting situation is that of a space curve
$X \subset \C \PP^3$. Our computations for space curves, both here and in \cite{RS}, were
performed in {\tt Macaulay2} \cite{M2},  but, even with ad hoc tricks, they turned out
to be more difficult than we had expected when we first embarked on our project.

Here is an illustration of the issue.
Let $X$ be the smooth sextic curve in $\C \PP^3$ defined by
$$ \bigl \langle \,x^2+y^2+z^2+w^2 , \, x y z-w^3 \, \bigr\rangle.
$$
The following lines of {\tt Macaulay2} code find the 
surface $X^*$ in $(\C \PP^3)^\vee$ that is dual to $X$:
\begin{verbatim}
S = QQ[x,y,z,w,X,Y,Z,W]; 
d = 4; pairing = first sum(d,i->(gens S)_i*(gens S)_{i+d});
makedual = I -> (e = codim I; J = 
saturate(I + minors(e+1,submatrix(jacobian(I+ideal(pairing)),{0..d-1},)),
minors(e,submatrix(jacobian(I),{0..d-1},)));eliminate((gens S)_{0..d-1},J))
makedual  ideal( x^2+y^2+z^2+w^2, x*y*z-w^3 );
\end{verbatim}
This program runs for a few minutes and outputs a polynomial
of degree $18$ with $318$ terms:
$$
729 x^{14} y^4+3861 x^{12} y^6+7954 x^{10} y^8+7954 x^8 y^{10}+3861 x^6 y^{12}+729 x^4 y^{14}+1458 x^{14} y^2 z^2 + \cdots
$$
Projective duality tends to produce large equations, even on modestly sized input,
and symbolic programs, like our little {\tt Macaulay2} fragment above, will often fail to terminate.

One promising alternative line of attack is offered by numerical algebraic geometry \cite{BHSW}.
Preliminary  experiments by Jonathan Hauenstein demonstrate that the software {\tt Bertini}
can perform the transformations $X \mapsto X^*$
and $X \mapsto (X^{[k]})^*$ in a purely numerical manner. 

Convex algebraic geometry requires the development of new
specialized software tools, both symbolic and numeric, and integrated
with optimization method. The advent of such new tools will
make our formula (\ref{mainformula}) more practical for non-linear convex hull
computations.

% \vfill \eject
 %\bigskip
\medskip

\noindent
{\bf Acknowledgments.}
This project started at the Banff International Research Station (BIRS)
during the workshop {\em Convex Algebraic Geometry} (February 14-18, 2010).
We are grateful to BIRS. Angelica Cueto and Herwig Hauser kindly allowed us to 
use their respective Figures 1 and 2.
We thank Roland Abuaf for his careful reading of the first version of this paper.
Bernd Sturmfels was supported in part by NSF grant DMS-0757207.

\bigskip

\end{document}